\documentclass[reqno]{amsart}

\usepackage{amsmath, amsthm, amscd, amsfonts, amssymb, graphicx, color,mathrsfs,mathtools,enumerate}
\usepackage{lineno}
\usepackage{dsfont}
\usepackage[bookmarksnumbered, colorlinks, plainpages]{hyperref}
\usepackage{bm}

\newtheorem{thm}{Theorem}[section]
\newtheorem{lemma}[thm]{Lemma}
\newtheorem{prop}[thm]{Proposition}
\newtheorem{coro}[thm]{Corollary}

\newtheorem{hyp1}[thm]{Hypotheses}

\newtheorem{rmk}[thm]{Remark}
\newtheorem{defn}[thm]{Definition}

\newcommand{\N}{\mathbb N}
\newcommand{\R}{\mathbb R}
\newcommand{\X}{\mathcal{X}}
\newcommand{\K}{\mathcal{K}}

\newcommand{\Id}{{\operatorname{I}}}
\newcommand{\norm}[1]{{\left\|#1\right\|}}
\newcommand{\scal}[2]{{\left\langle #1,#2\right\rangle}}

\newcommand{\E}{{\mathbb E}}
\newcommand{\abs}[1]{{\left|#1\right|}}

\allowdisplaybreaks

\title{Weak uniqueness for stochastic partial differential equations in Hilbert spaces}
\author{Davide Addona, Davide Augusto Bignamini}

\address{D.Addona: Dipartimento di Scienze Matematiche, Fisiche e Informatiche, Universit\`a degli Studi di Parma, Parco Area delle Scienze, 53/A, 43124 Parma, Italy PR}
\email{\textcolor[rgb]{0.00,0.00,0.84}{davide.addona@unipr.it}}

\address{D.A. Bignamini: Dipartimento di Scienza e Alta Tecnologia (DISAT), Universit\`a degli Studi dell'In\-su\-bria, Via Valleggio 11, 22100 Como, Italy}
\email{\textcolor[rgb]{0.00,0.00,0.84}{da.bignamini@uninsubria.it}}

\keywords{weak uniqueness; regularization by noise; unbounded continuous drift; directional gradients; stochastic damped wave equation; stochastic heat equation \\$ $\\ 
D.Addona: Dipartimento di Scienze Matematiche, Fisiche e Informatiche, Universit\`a degli Studi di Parma, Parco Area delle Scienze, 53/A, 43124 Parma, Italy PR, davide.addona@unipr.it\\
D.A. Bignamini: Dipartimento di Scienza e Alta Tecnologia (DISAT), Universit\`a degli Studi dell'In\-su\-bria, Via Valleggio 11, 22100 Como, Italy, da.bignamini@uninsubria.it}
\thanks{The authors are members of GNAMPA (Gruppo Nazionale per l’Analisi Matematica, la Probabilit\`a e le loro Applicazioni) of the Istituto Nazionale di Alta Matematica (INdAM)}
\subjclass[2020]{60H15,60H50}

\begin{document}
\begin{abstract}
Let $U,H$ be two separable Hilbert spaces. The main goal of this paper is to study the weak uniqueness of the Stochastic Differential Equation evolving in $H$
\begin{align*}
dX(t)=AX(t)dt+\mathcal{V}B(X(t))dt+GdW(t), \quad t>0, \quad X(0)=x \in H,    
\end{align*}
where $\{W(t)\}_{t\geq 0}$ is a $U$-cylindrical Wiener process, $A:D(A)\subseteq H\to H$ is the infinitesimal generator of a strongly continuous semigroup, $\mathcal{V},G:U\rightarrow H$ are linear bounded operators and $B:H\rightarrow U$ is a uniformly continuous function. The abstract result in this paper gives the weak uniqueness for large classes of heat and damped equations in any dimension without any H\"older continuity assumption on $B$.
\end{abstract}

\maketitle

\section*{Introduction}
Let $d\in\N$. The main goal of this paper is to investigate the weak uniqueness of the following class of stochastic damped equations 
\begin{align}\label{damped_stoc_equationI}
\left\{
\begin{array}{ll}
\displaystyle \frac{\partial^2y}{\partial t^2}(t)
= -\Lambda y(t)-\rho\Lambda^{\alpha}\left(\frac{\partial y}{\partial t}(t)\right)+\Lambda^{-\beta}B\left(y(t),\frac{\partial y}{\partial t}(t)\right)+\Lambda^{-\gamma}d{W}(t), & t\in[0,T], \vspace{1mm} \\
y(0)=y_0\in L^2([0,2\pi]^d), & \vspace{1mm} \\
\displaystyle \frac{\partial y}{\partial t}(0)=y_1\in L^2([0,2\pi]^d),
\end{array}
\right.
\end{align}
where $\Lambda=(-\Delta)^k$, $\Delta$ is the realization of the Laplace operator with Dirichlet boundary conditions, $k=1,2$, $\alpha\in(0,1)$, $\rho,\beta,\gamma\geq0$, $B$ is a uniformly continuous function and $\{W(t)\}_{t\geq 0}$ is a $L^2([0,2\pi]^d)$-cylindrical Wiener process. Pathwise uniqueness for SPDEs such as \eqref{damped_stoc_equationI} has already been studied in the literature (see, for instance, \cite{AddBig, AddBig2, AddMasPri23}, we also refer to \cite{MasPri2017,MasPri2024} for the case $\alpha=\rho=0$). However, to the best of our knowledge, there is not any result available on weak uniqueness for such types of SPDEs. In \cite{Han2024}, weak uniqueness for \eqref{damped_stoc_equationI} is investigated in the case $d=1$, $\alpha=0$ and the diffusion coefficient $\Lambda^{-\gamma}$ is replaced by a function $G$ that depends on the solution. We note that the assumption \cite[H5]{Han2024} on $G$ does not allow for considering $G \equiv \Lambda^{-\gamma}$ as in this paper. Also in \cite{Cho-Gol1995,Kun2013} the weak well-posedness for abstract SPDEs is studied but, in both papers, the controllability condition assumed on the coefficients excludes SPDEs such as \eqref{damped_stoc_equationI} (see Remarks \ref{diff-kunze} and \ref{rmk-damped}).

Here, we aim to prove weak uniqueness for SPDEs such as \eqref{damped_stoc_equationI} under weaker conditions than those assumed for pathwise uniqueness in \cite{AddBig, AddBig2, AddMasPri23}. In particular, under suitable assumptions on $\alpha, \gamma, \beta$, we will prove that for every $d \in \mathbb{N}$, weak uniqueness holds for \eqref{damped_stoc_equationI} without any additional Hölder continuity condition on $B$.
We will operate within an abstract framework and, by taking advantage of the finite-dimensional approximations technique presented in \cite{AddBig}, we notably simplify, in this setting, the method presented in \cite{BerOrrSca2024} for the stochastic heat equation. Inspired by \cite{MasPri2017,MasPri2024}, we will relax the controllability assumptions on the coefficient of \eqref{damped_stoc_equationI} exploiting the regularizing properties of Ornstein-Uhlenbeck type semigroups along suitable directions (see for instance \cite{BigFer2023,BigFerForZan2023,CanDap1996}).
We stress that deterministic linear damped operator, corresponding to the linear part of \eqref{damped_stoc_equationI}, has been widely investigated, also in more general situations ( see \cite{CheRus1981,CheTri1988,CheTri1989,CheTri1990-1,CheTri1990-2, LasTri1998,Tri2003}). We stress that, in \cite{Tri2003}, it has been proved the optimal rate of explosion at $0$ of the control problem associated to such an operator. In the present paper, due to the presence of the operator $\Lambda^{-\beta}$ in \eqref{damped_stoc_equationI}, we will exploit the spectral approach presented in the aforementioned paper to establish similar estimates along appropriate directions.

Finally, the abstract result in this paper improves those in \cite{Cho-Gol1995,Kun2013} for a class of stochastic heat equations.  We refer to \cite{BerOrrSca2024,Pri2021} for weak uniqueness results for a class of stochastic heat equations, in the case when the nonlinear part of the drift takes values in space of distributions.

\section{Notation}
\label{sec:notation}
 Let $\K$ be a Banach space endowed with the norm $\norm{\cdot}_\K$. We denote by $\mathcal{B}(\K)$ the Borel $\sigma$-algebra associated to the norm topology in $\K$. We denote by $\mathbb{E}$ the expectation with respect to a probability measure $\mathbb{P}$.  

Let $\K_1$ and $\K_2$ be two real separable Banach spaces equipped with the norms $\norm{\cdot}_{\K_1}$ and $\norm{\cdot}_{\K_2}$, respectively.
We denote by $\Id_{\K_1}$ the identity operator on $\K_1$.
We denote by $\mathcal{L}(\K_1;\K_2)$ the set of bounded linear operator from $\K_1$ to $\K_2$, if $\K_1=\K_2$ then we write $\mathcal{L}(\K_1)$.

We denote by $B_b(\K_1;\K_2)$ the set of bounded and Borel measurable functions from $\K_1$ into $\K_2$. For every $f\in B_b(\K_1;\K_2)$ we set $\norm{f}_\infty=\sup_{x\in\K_1}\norm{f(x)}_{\K_2}.$
We denote by $UC(\K_1;\K_2)$ the space of uniformly continuous functions from $\K_1$ into $\K_2$ and by $UC_b(\K_1;\K_2)$ its subspace of bounded and  uniformly continuous functions.  We denote by $UC^k_b(\K_1;\K_2)$ the space of $k$-times Fr\'echet differentiable functions bounded, uniformly continuous with bounded and uniformly continuous derivatives up to order $k$. We denote by $D^i$ the $i$-order Fr\'echet derivative, with $i\in\N$. We endow $B_b(\K_1;\K_2)$, $UC_b(\K_1;\K_2)$ and  $UC^k_b(\K_1;\K_2)$ with the natural norms. We set $UC^\infty_b(\K_1;\K_2)=\cap_{k\in\N}UC^k(\K_1;\K_2)$ and, if $\K_2=\R$, then we simply write $B_b(\K_1)$, $UC_b(\K_1)$,  $UC^k_b(\K_1)$ and $UC^\infty_b(\K_1)$. 

Let $\X$ be a separable Hilbert space. If $f\in UC^2_b(\X)$, then we denote by $\nabla f (x)$ and $\nabla^2f(x)$ the Fr\'echet gradient and Hessian operator at $x\in\X$, respectively.

\section{Abstract framework and main result}
Let $U,H$ be two separable Hilbert spaces. We consider an SDE which evolves in the Hilbert space $H$ of the form
\begin{align}
\label{SDE}
dX(t)=AX(t)dt+\mathcal{V}B(X(t))dt+GdW(t), \quad t>0, \quad X(0)=x \in H,    
\end{align}
whose coefficients satisfy the following assumptions.


\begin{hyp1}\label{hyp:standard}
Assume that the following conditions hold true.
\begin{enumerate}[\rm(I)]
\item $\{W(t)\}_{t\geq 0}$ is a $U$-cylindrical Wiener process.

\item $A:{\rm Dom}(A)\subseteq H\rightarrow H$
is the infinitesimal generator of a compact strongly continuous semigroup. $\mathcal{V}, G\in\mathcal L(U;H)$ are linear bounded operators and $B\in UC(H;U)$.

\item There exists $\eta\in (0,1)$ such that for every $t>0$ we have
\begin{equation}\label{tra-con}
\int^t_0\frac{1}{s^\eta}{\rm Trace}_H\left[e^{sA}GG^*e^{sA^*}\right]ds<\infty.
\end{equation}
\end{enumerate}
\end{hyp1}

\begin{defn}\label{def-sol}
Let $x\in H$.

\begin{itemize}
\item[(Existence)] A weak mild solution to \eqref{SDE} is a pair $(X,W)$ where $W$ is a $U$-cylindrical Wiener process defined on a complete filtered probability space 
$(\Omega,\mathcal{F},\{\mathcal{F}_t\}_{t\geq0},\mathbb{P})$ and 
$X$ is a $H$-valued $\{\mathcal{F}_t\}_{t\geq0}$-adapted process having continuous trajectories $P$-a.s. and such that, for every $t\geq0$, 
\begin{align}\label{mild}
X(t)=e^{tA}x+\int_0^te^{(t-s)A}\mathcal{V}B(X(s))\, d s+W_A(t),
\qquad \mathbb{P}\text{-a.s.},   
\end{align}
where $W_A$ is the stochastic convolution process given by 
\begin{equation}\label{W_A}
    W_A(t):=\int^t_0e^{(t-s)A}G\, dW(s)\,, \quad t\geq 0.
\end{equation}
\item[(Uniqueness)] We say that weak uniqueness holds for \eqref{SDE} if, whenever $(X,W)$ and $(Y,W)$ are two mild solutions to \eqref{SDE}, then $X$ and $Y$ have the same law on $C(\R_+;H)$, namely that for every bounded $\psi:C(\R_+;H)\to \R$ we have
\[
\E\left[\psi(X)\right] = \E\left[\psi(Y)\right].
\]
\end{itemize}
\end{defn}

We recall the following result about the existence of a weak solutions.

\begin{prop}{\cite[Proposition 3]{Cho-Gol1995}}
Assume that Hypotheses \ref{hyp:standard} hold and $B$ has linear growth. Then \eqref{SDE} admits a weak mild solution.
\end{prop}
In this paper, we are going to prove that weak uniqueness holds true for \eqref{SDE} under the following additional conditions.

\begin{hyp1}\label{hyp:uniqueness}
Assume Hypotheses \ref{hyp:standard} and that the following conditions are satisfied.
\begin{enumerate}[\rm(i)]
\item There exists a sequence of finite-dimensional subspaces $\{H_n\}_{n\in\N}\subseteq H$ such that \\ $H=\overline{\cup_{n\in\N}H_n}$, $H_0:=\{0\}$ and for every $n\in\N$ we have
\begin{align*}
& H_{n-1}\subseteq H_{n}, \qquad 
H_{n}\subseteq {\rm Dom}(A)\cap {\rm Dom}(A^*) ,\qquad A(H_{n}\backslash H_{n-1})\subseteq \left(H_{n}\backslash H_{n-1}\right)\cup \{0\}.
\end{align*}
\item For every $t>0$ we have
\begin{align}
&e^{tA}(H)\subseteq Q^{\frac12}_t(H),\qquad\qquad \qquad\qquad\quad Q_{t}:=\int^t_0e^{sA}GG^*e^{sA^*}ds,
   \label{contron}\\
&\int^{t}_0 \|\Gamma_s\mathcal{V}\|_{\mathcal{L}(U;H)}ds<\infty, \qquad\qquad\qquad\Gamma_s:=Q^{-\frac12}_se^{sA}, \ s>0.
\label{supercontron}
\end{align}
\end{enumerate}

\end{hyp1}
\begin{rmk}\label{diff-kunze}
We underline that, in many interesting cases, condition \eqref{supercontron} is weaker than the condition 
\begin{equation}\label{K-contro}
\int^t_0 \|\Gamma_s\|_{\mathcal{L}(H;H)}ds<\infty, \qquad t>0,
\end{equation}
assumed in  \cite{Cho-Gol1995,Kun2013} to guarantee weak uniqueness. For example in the case of stochastic damped equations presented in Subsection \ref{sub:stoch_damp_eq}, condition \eqref{K-contro} is never satisfied (see Remark \ref{rmk-damped}). Moreover, in the case of stochastic heat equations in $H=L^2([0,2\pi]^d)$, assuming \eqref{tra-con} and \eqref{K-contro} implies $d<4$ (see Remark \ref{rmk-heat}). On the contrary, assumption \eqref{supercontron} in the present paper allows us to avoid these issues. Indeed, we will establish weak well-posedness for a large class of stochastic damped equations (see Theorems \ref{Wave} and \ref{EB}), as well as for a class of heat equations in every dimension (see Theorem \ref{Heat-weak}). We emphasize that to leverage condition \eqref{supercontron} instead of \eqref{K-contro}, we will consider a notion of differentiability along the range of $\mathcal{V}$, see Definition \ref{C-D_diff}.
\end{rmk}

Now we state the main result of this paper.
\begin{thm}\label{WeakLoc}
Assume that Hypotheses \ref{hyp:uniqueness}  are fulfilled. Then weak uniqueness holds  \eqref{SDE}.
\end{thm}
\begin{rmk}
We will prove Theorem \ref{WeakLoc} in the case $B\in UC_b(H;U)$ (see Theorem \ref{Weak1}). The procedure presented in \cite[Section 6]{Pri2021} yields Theorem \ref{WeakLoc} in its full generality.
\end{rmk}

\begin{rmk}
We note that, whereas the assumptions for the pathwise uniqueness results in \cite{AddBig,AddBig2,AddMasPri23} impose an additional H\"older condition on $B$, the present result on weak uniqueness does not. This is due to the fact that, for weak uniqueness, a uniform estimate for the second-order derivative of the solution to \eqref{eq:n-staz} is not required (see \cite{BerOrrSca2024, Pri2021}).
\end{rmk}

\section{Proof of the main result}

\subsection{Ornstein-Uhlenbeck semigroup}

We introduce the Ornstein-Uhlenbeck semigroup $\{R(t)\}_{t\geq 0}$ given by
\begin{equation}\label{OUS}
(R(t)\varphi)(x)=\int_{H}\varphi(e^{tA}x+y)\mu_t(dy),\quad t>0,\;\varphi\in B_b(H),\; x\in H,
\end{equation}
where $\mu_t$ is the Gaussian measure on $\mathcal{B}(H)$ with mean $0$ and covariance operator $Q_t$ given by \eqref{contron}. First of all we recall that for every $\varphi\in B_b(H)$ and $t>0$ we have
\begin{equation}\label{zero-zero}
\norm{R(t)\varphi}_{\infty}\leq \norm{\varphi}_{\infty}.
\end{equation}
Under Hypotheses \ref{hyp:standard} the semigroup $\{R(t)\}_{t\geq 0}$ verifies some regularity properties that we state in the subsequent propositions.
\begin{prop}{\cite[Theorem 6.2.2 and Proposition 6.2.9]{Dap-Zab2002}}
\label{stimeclassiche}
Assume that Hypotheses \ref{hyp:standard} hold. For every $t>0$ we have
\[
R(t)(B_b(H))\subseteq UC_b^{\infty}(H).
\]
In particular, for every $\varphi\in B_b(H)$, $t>0$ and $x,h\in H$ we have
\begin{align}
&DR(t)\varphi(x)h=\int_{H}\scal{\Gamma_{t}h}{Q_{t}^{-\frac12}y}_H\varphi(e^{tA}x+y)\mu_{t}(dy),\label{D10}\\
&|DR(t)\varphi(x)h|\leq \|\Gamma_th\|_H\|\varphi\|_{\infty}.\label{stima10}
\end{align}
\end{prop}

To exploit condition \eqref{supercontron} instead of \eqref{K-contro} we need to introduce a notion of regularity along suitable directions, see for instance \cite{BigFer2023,BigFerForZan2023,CanDap1996}.

\begin{defn}\label{C-D_diff}
We say that a function $\varphi\in UC_b(H)$ is $\mathcal{V}$-differentiable at $x\in H$ if there exists $\ell_x\in U$ such that for every $v\in U$ it holds 
\begin{equation}\label{C-D_1}
\lim_{s\to 0}\abs{\frac{f(x+s\mathcal{V}v)-f(x)}{s}-\langle \ell_x,v\rangle_U}=0.
\end{equation}
We set $\nabla_{\mathcal{V}}f(x):=\ell_x$. We say that a function is $\mathcal{V}$-differentiable if it is $\mathcal{V}$-differentiable at any $x\in H$. We denote by $C^1_{b,\mathcal{V}}(H)$ the subspace of $UC_b(H)$ of $\mathcal{V}$-differentiable functions $f$ such that $\nabla_{\mathcal{V}}f\in UC_b(H;U)$.
\end{defn}
We note that $C^{1}_{b,\mathcal{V}}(H)$ is a Banach space endowed with the norm
\[
\norm{f}_{C^{1}_{b,\mathcal{V}}(H)}:=\norm{f}_{\infty}+\norm{\nabla_\mathcal{V}f}_{\infty}.
\]
In the next proposition we underline the relation between $\mathcal{V}$- and G\^ateaux differentiability.

\begin{prop}
\label{gat-Rdiff}
Assume that Hypotheses \ref{hyp:standard} hold true. If $\varphi\in UC_b(H)$ is G\^ateaux differentiable, then $\varphi$ is $\mathcal{V}$-differentiable and $\nabla_{\mathcal{V}}\varphi(x)=\mathcal{V}^*\nabla\varphi(x)$ for every $x\in H$.
In particular, $UC^1_b(H)\subseteq C^1_{b,\mathcal{V}}(H)$.
\end{prop}
\begin{proof}
Let $\varphi:H\to\R$ be a G\^ateaux differentiable function and let $x\in H$ and $v\in U$. Since $\mathcal{V}v\in H$, by the definition of G\^ateaux differentiability we have
\begin{equation*}
\lim_{s\to 0}\abs{\frac{\varphi(x+s\mathcal{V}v)-\varphi(x)}{s}-\langle \mathcal{V}^*\nabla\varphi(x),v\rangle_U}=\lim_{s\to 0}\abs{\frac{\varphi(x+s\mathcal{V}v)-\varphi(x)}{s}-\langle \nabla\varphi(x),\mathcal{V}v\rangle_H}=0.
\end{equation*}
Hence, $\mathcal{V}^*\nabla\varphi(x)$ verifies Definition \ref{C-D_diff}, and so $\varphi$ is $\mathcal{V}$-differentiable with $\nabla_{\mathcal{V}}\varphi=\mathcal{V}^*\nabla\varphi$.
\end{proof}
By \cite[Proposition 3.3]{LunRoc2021}, for every $\varphi\in C^{1}_{b,\mathcal{V}}(H)$ and $t>0$ we have
\begin{equation}\label{theta-theta}
\norm{R(t)\varphi}_{C^{1}_{b,\mathcal{V}}(H)}\leq \left(1+\norm{e^{tA^*}}_{\mathcal{L}(H;H)}\right)\norm{\varphi}_{C^{1}_{b,\mathcal{V}}(H)}.
\end{equation}
By \eqref{zero-zero} and Proposition \ref{stimeclassiche} (with $h=\mathcal{V}v$ and $v\in U$) we immediately deduce the following result.
\begin{prop}\label{Vschauder}
Assume that Hypotheses \ref{hyp:standard} hold true. For every $\varphi\in UC_b(H)$ and $t>0$ we have
\[
\norm{R(t)\varphi}_{C^{1}_{b,\mathcal{V}}(H)}\leq \left(1+\norm{e^{tA^*}}_{\mathcal{L}(H)}\norm{\Gamma_t\mathcal{V}}_{\mathcal{L}(U,H)}\right)\norm{\varphi}_{\infty}.
\]
\end{prop}

Let $g\in UC_b(H)$ and $\lambda>0$. We consider the following integral equation
\begin{equation}\label{eq:staz}
u(x)=\int_0^{\infty}e^{-\lambda t}R(t)\left[\scal{\nabla_{\mathcal{V}}u}{B}_U+g\right](x)dt,\qquad x\in H.
\end{equation}
\begin{prop}\label{Hstaz}
Assume that Hypotheses \ref{hyp:uniqueness}  hold and $B\in UC_b(H;U)$. There exists $\lambda_0>0$ such that for every $\lambda>\lambda_0$ and $g\in UC_{b}(H)$ equation \eqref{eq:staz} admits a unique solution $u_\lambda\in C^{1}_{b,\mathcal{V}}(H)$ and
\begin{equation}\label{Vu-stima}
\norm{u_\lambda}_{C^{1}_{b,\mathcal{V}}(H)}\leq \norm{g}_{\infty}.
\end{equation}
\end{prop}
\begin{proof}
For every $\lambda>0$ we consider the operator $V_\lambda$ defined, for every $u\in C^{1}_{b,\mathcal{V}}(H)$, by
\[
V_\lambda(u)=\int_0^{\infty}e^{-\lambda t}R(t)\left[\scal{\nabla_{\mathcal{V}}u}{B}_U+g\right](x)dt.
\]
To conclude the proof it is sufficient to prove that $V_\lambda$ is a contraction in $C^{1}_{b,\mathcal{V}}(H)$ for every $\lambda>\lambda_0$, where $\lambda_0$ is a suitable positive constant. By \eqref{theta-theta} and Proposition \ref{Vschauder}, for every $u\in C^{1}_{b,\mathcal{V}}(H)$ and $\lambda>0$ we have
\begin{align}
\norm{V_\lambda(u)}_{C^{1}_{b,\mathcal{V}}(H)}&\leq \int_0^{1}e^{-\lambda t}\norm{R(t)\left[\scal{\nabla_{\mathcal{V}}u}{B}_U+g\right]}_{C^{1}_{b,\mathcal{V}}(H)}dt\notag\\
&+\int_1^{\infty}e^{-\lambda t}\norm{R(t-1)R(1)\left[\scal{\nabla_{\mathcal{V}}u}{B}_U+g\right]}_{C^{1}_{b,\mathcal{V}}(H)}dt\notag\\
&\leq \left(\norm{\nabla_{\mathcal{V}}u}_{ \infty}\norm{B}_{\infty}+\norm{g}_{\infty}\right)\int_0^{1}e^{-\lambda t}\left(1+\norm{e^{tA^*}}_{\mathcal{L}(H)}\norm{\Gamma_t\mathcal{V}}_{\mathcal{L}(U;H)}\right)dt,\notag\\
&+\int_1^{\infty}e^{-\lambda t}\left(1+\norm{e^{tA^*}}_{\mathcal{L}(H)}\right)\norm{R(1)\left[\scal{\nabla_{\mathcal{V}}u}{B}_U+g\right](x)}_{C^{1}_{b,\mathcal{V}}(H)}dt\notag\\
&\leq \left(\norm{\nabla_{\mathcal{V}}u}_{ \infty}\norm{B}_{\infty}+\norm{g}_{\infty}\right)(c_1+c_2),\label{pre-contraz}
\end{align}
where
\begin{align}
&c_1:= \int_0^{1}e^{-\lambda t}\left(1+\norm{e^{tA^*}}_{\mathcal{L}(H)}\norm{\Gamma_t\mathcal{V}}_{\mathcal{L}(U;H)}\right)dt\label{c1},\\
&c_2:=\left(1+\norm{e^{A^*}}_{\mathcal{L}(H)}\norm{\Gamma_1\mathcal{V}}_{\mathcal{L}(U;H)}\right)\int_0^{\infty}e^{-\lambda t}\left(1+\norm{e^{tA^*}}_{\mathcal{L}(H)}\right)dt\label{c2}.
\end{align}
Since $\{e^{tA}\}_{t\geq 0}$ is strongly continuous on $H$, by \eqref{supercontron} the constants $c_1$ and $c_2$ are finite so $V_\lambda$ is bounded from $C^{1}_{b,\mathcal{V}}(H)$ into itself. By \eqref{pre-contraz} and Proposition \ref{Vschauder} for every $u,z\in C^{1}_{b,\mathcal{V}}(H)$ and $\lambda>0$ we have
\begin{align*}
\norm{V_\lambda(u)-V_\lambda(z)}_{ C^{1}_{b,\mathcal{V}}(H)}&\leq \int_0^{1}e^{-\lambda t}\norm{R(t)\left[\scal{\nabla_{\mathcal{V}}u-\nabla_{\mathcal{V}}z}{B}_U\right](x)}_{C^{1}_{b,\mathcal{V}}(H)}dt\\
&+\int_1^{\infty}e^{-\lambda t}\norm{R(t-1)R(1)\left[\scal{\nabla_{\mathcal{V}}u-\nabla_{\mathcal{V}}z}{B}_U\right](x)}_{C^{1}_{b,\mathcal{V}}(H)}dt\\
&\leq \norm{u-z}_{C^{1+\theta}_{b,\mathcal{V}}(H)}\norm{B}_{\infty}(c_1+c_2),
\end{align*}
where $c_1$ and $c_2$ are given by \eqref{c1} and \eqref{c2}, respectively. Hence by condition \eqref{supercontron} there exists $\lambda_0>0$ such that for every $\lambda>\lambda_0$ we have
\begin{equation}\label{contraz}
c_1+c_2<\frac{1}{2},\qquad\norm{B}_{\infty}(c_1+c_2)<\frac{1}{2}.
\end{equation}
This implies that $V_\lambda$ is a contraction on $ C^{1}_{b,\mathcal{V}}(H)$. Finally, \eqref{Vu-stima} follows from \eqref{pre-contraz} when \eqref{contraz} holds.
\end{proof}

\subsection{Finite-dimensional approximations}
Assume that Hypotheses \ref{hyp:uniqueness}  hold true and that $B\in UC_b(H;U)$. For every $n\in\N$, we set
\begin{align}\label{coefficienti-n}
A_n:=AP_n,\quad \mathcal{V}_n:=P_n \mathcal{V},\quad G_n=P_nG, \quad B_n=B\circ P_n,
\end{align}
where $P_n$ is orthogonal projection on $H_n$.
By \cite[Lemma 4.1]{AddBig} for every $n\in\N$ we have
\begin{align}
& e^{tA}_{|H_n}=e^{tA_n} \qquad \forall t\geq0,\label{P5}\\
&\norm{B_n}_{\infty}\leq \norm{B}_{\infty}.
\label{holdrBn}
\end{align}
Fix $x\in H$ and let $(X,W)$ be a weak mild solution to \eqref{SDE}.
For every $n\in\N$ we introduce the process $X_n=\{X_n(t,x)\}_{t\geq 0}$ defined, for every $t\geq 0$, by
\begin{align}
\label{mild_sol_n}
X_n(t,x)=e^{tA_n}P_nx+\int_0^t e^{(t-s)A_n}\mathcal{V}_n B_n(X(s,x))ds+W_{A,n}(t),\qquad \mathbb P{\rm -a.s.,}
\end{align}
where 
\begin{equation}\label{csn}
W_{A,n}(t)=\int_0^te^{tA_n}G_ndW(s), \qquad t\geq 0.
\end{equation}
In particular, $X_n$ fulfills in the standard It\^o sense
\begin{align}
\label{mild_sol_n_var}
& dX_n(t,x)=A_nX_n(t,x)dt+\mathcal{V}_nB_n(X(t,x))dt+G_ndW(t), \quad t\in[0,T], \qquad  X_n(0,x)=P_nx. \end{align}

\begin{lemma}{\cite[Lemma 4.2]{AddBig}}
Assume that Hypotheses \ref{hyp:uniqueness} hold and $B\in UC_b(H;U)$. For every fixed $T>0$ and $x\in H$, we have
\begin{equation}\label{supE}
\lim_{n \to \infty}\sup_{t\in [0,T]}\mathbb{E}\left[\|X_{n}(t,x)-X(t,x)\|_H^2\right]dt=0.
\end{equation} 
\end{lemma}
\begin{proof}
The proof is analogous to that of \cite[Lemma 4.2]{AddBig}, which is still valid if $B\in UC_b(H;U)$ and it is not assumed to be H\"older continuous. 
\end{proof}

For every $n\in\N$, we define the finite-dimensional Ornstein-Uhlenbeck semigroup $\{R_n(t)\}_{t\geq 0}$ as
\[
(R_n(t)\varphi)(x)=\int_{H_n}\varphi(e^{tA_n}x+y)\mu_{t,n}(dy),\quad t>0,\;\varphi\in B_b(H_n),\; x\in H_n,
\]
where $\mu_{t,n}$ is the Gaussian measure on $\mathcal{B}(H_n)$ with mean $0$ and covariance operator 
\[
Q_{t,n}:=\int_0^te^{sA_n}G_nG_n^*e^{sA_n^*}ds, \qquad t>0.
\]
For every $t> 0$ and $n\in\N$ we set
\[
\Gamma_{t,n}:=Q_{t,n}^{-1/2}e^{tA_n}.
\]
By \cite[Proposition 4.4]{AddBig}, conditions \eqref{contron} and \eqref{supercontron} imply that for every $t>0$ we have
\begin{align}
&e^{tA_n}(H_n)\subseteq Q_{t,n}^{\frac12}(H_n),\qquad n\in\N,\label{Qtn-iniettivo}\\
&\sup_{n\in\N}\int_0^{t}\norm{\Gamma_{s,n}\mathcal{V}_n}_{\mathcal{L}(U;H_n)}ds<\infty\label{Cesistenza}.
\end{align}    
Let $n\in\N$, $g\in UC_b(H_n)$ and $\lambda>0$. We consider the following integral equation
\begin{equation}\label{eq:staz-n}
u(x)=\int_0^{\infty}e^{-\lambda t}R_n(t)\left[\scal{\nabla_{\mathcal{V}_n}u}{B_n}_U+g\right](x)dt,\qquad x\in H_n.
\end{equation}

\begin{prop}\label{Sol-n-int}
Assume that Hypotheses \ref{hyp:uniqueness}  hold true and that $B\in UC_b(H;U)$. There exists $\lambda_0>0$ such that for every $n\in\N$, $\lambda>\lambda_0$ and $g\in UC_b(H_n)$ the equation \eqref{eq:staz-n} has unique solution $u\in C^{1}_{b,\mathcal{V}_n}(H_n)$ and
   \begin{equation*}
\norm{u_{\lambda,n}}_{C^{1}_{b,\mathcal{V}_n}(H_n)}\leq \norm{g}_{\infty}.
  \end{equation*}
\end{prop}
\begin{proof}
The statements follow by \eqref{P5}, \eqref{holdrBn}, \eqref{Qtn-iniettivo}, \eqref{Cesistenza} and applying Proposition \ref{Hstaz} with $H$, $R(t)$, $\mathcal{V}$ and $B$ replaced by $H_n$, $R_n(t)$, $\mathcal{V}_n$ and $B_n$, respectively. We underline that $\lambda_0$ is independent of $n\in\N$ thanks to \eqref{Cesistenza}.
\end{proof}
\begin{prop}\label{Sol-nh}
Assume that Hypotheses \ref{hyp:uniqueness}  hold and $B\in UC_b(H;U)$. Let $\lambda_0$ be the constant given by Proposition \ref{Sol-n-int}. Fix $n\in\N$, $g\in UC_b(H_n)$ and $\lambda>\lambda_0$. There exists a sequence $(f_h)_{h\in\N}\subseteq UC^{2}_b(H_n)$ such that for every $h\in\N$, the following elliptic equation 
    \begin{equation}\label{eq:n-staz}
\lambda u(x) -\frac{1}{2}{\rm Trace}\left[G_nG_n^*\nabla^2u(x)\right]
  +\scal{A_nx}{\nabla u(x)}
    =f_h(x)\,, \quad x\in H_n,
 \end{equation}
     has a unique strict solution $u_{\lambda,n,h}\in UC^{2}_b(H_n)$ given by
\[
u_{\lambda,n,h}(x)=\int_0^{\infty}e^{-\lambda t}R_n(t)f_h(x)dt;
\]
Moreover
\begin{align}
&\sup_{h\in\N}\|f_h\|_{\infty}<\infty, \qquad \lim_{h\rightarrow\infty}\norm{f_h-\scal{\nabla_{\mathcal{V}_n}u_{\lambda,n}}{B_n}_U-g_n}_{\infty}=0,\label{conv-g}\\
&\lim_{h\rightarrow \infty}\norm{u_{\lambda,n,h}-u_{\lambda,n}}_{C^1_{b,\mathcal{V}_n}(H_n)}=0,\label{app-sol-par}
\end{align}
where $u_{\lambda,n}$ is the unique solution to \eqref{eq:staz-n} given by Proposition \ref{Sol-n-int}.
\end{prop}

\begin{proof}
Let $\lambda_0$ be the constant given by Proposition \ref{Sol-n-int}, let $\lambda>\lambda_0$, $n\in\N$, $g\in UC_b(H_n)$  and $u_{\lambda,n}$ be the unique solution to \eqref{eq:staz-n} given by Proposition \ref{Sol-n-int}. By Proposition \eqref{Sol-n-int} we know that
 \[
\scal{\nabla_{\mathcal{V}_n}u_{\lambda,n}}{B_n}_U+g_n\in UC_b(H_n),
 \]
so by standard approximation arguments ($f_h$ should be defined by means of convolution with mollifiers), \eqref{conv-g} is verified. All the statements are quite classical, we refer to \cite{Lun1997} for a proof and a detailed discussion. The only difference is that in \cite{Lun1997} the convergence result \eqref{app-sol-par} is proved with respect to the $C_b$-norm, so here we prove a finer result in this particular case. 
By Proposition \ref{Vschauder}, and arguing as in \eqref{pre-contraz} there exists a constant $C>0$ such that 
\begin{align*}
\norm{u_{\lambda,n,h}(t,\cdot)-u_{\lambda,n,}(t,\cdot)}_{C^{1}_{b,\mathcal{V}_n}(H_n)}\leq C\norm{f_h-\scal{\nabla_{\mathcal{V}_n}u_{\lambda,n}}{B_n}_U-g_n}_{\infty}, \qquad h\in\N,
\end{align*}
so \eqref{conv-g} yields \eqref{app-sol-par}.
\end{proof}

\subsection{Weak uniqueness}
\begin{thm}\label{Weak1}
Assume that Hypotheses \ref{hyp:uniqueness}  hold true and $B\in UC_b(H;U)$. Then, weak uniqueness holds for \eqref{SDE}.
\end{thm}
\begin{proof}
Fix $x\in H$. Let $(X,W)$ and $(Y,W)$ be to weak mild solutions to \eqref{SDE} and let $X_n$ and  $Y_n$ be the finite-dimensional approximations of $X$ and $Y$ given by \eqref{mild_sol_n}, respectively. Let $g\in UC_b(H)$ and let $\lambda_0$ be the constant given by Proposition \ref{Sol-n-int}. For every $\lambda>\lambda_0$ and $n\in\N$, let $u_{\lambda,n}$ be the solution to \eqref{eq:staz-n} with $g$ replaced by $g_{|H_n}$. Let $\{f_h\}_{h\in\N}$ and $\{u_{\lambda,n,h}\}_{h\in \N}$ be the sequences given by Proposition \ref{Sol-nh}. Applying the It\^o formula to $e^{-\lambda t} u_{\lambda,n,h}(X_n)$ and exploiting \eqref{mild_sol_n_var} and \eqref{eq:n-staz}, we get
\begin{align*}
\mathbb{E}\left[e^{-\lambda t}u_{\lambda,n,h}(X_n(t))\right]&-u_{\lambda,n,h}(P_nx)+\int_0^te^{-\lambda s}\mathbb{E}\left[f_h(X_n(s))\right]ds\\
&=\int_0^te^{-\lambda s}\mathbb{E}\left[\scal{B_n(X(s))}{\nabla_{\mathcal{V}_n} u_{\lambda,n,h}(X_n(s))}_U\right]ds, \qquad t>0.
\end{align*}
Now letting $h\rightarrow \infty$ by \eqref{conv-g} and \eqref{app-sol-par} we obtain
\begin{align*}
\mathbb{E}\left[e^{-\lambda t}u_{\lambda,n}(X_n(t))\right]&-u_{\lambda,n}(P_nx)+\int_0^te^{-\lambda s}\mathbb{E}\left[g(X_n(s))\right]ds\\
&=\int_0^te^{-\lambda s}\mathbb{E}\left[\scal{B_n(X(s))-B_n(X_n(s))}{\nabla_{\mathcal{V}_n} u_{\lambda,n}(X_n(s))}_U\right]ds, \qquad t>0.
\end{align*}
In the same way, for every $n\in\N$, we obtain 
\begin{align*}
\mathbb{E}\left[e^{-\lambda t}u_{\lambda,n}(Y_n(t))\right]&-u_{\lambda,n}(P_nx)+\int_0^te^{-\lambda s}\mathbb{E}\left[g(Y_n(s))\right]ds\\
&=\int_0^te^{-\lambda s}\mathbb{E}\left[\scal{B_n(Y(s))-B_n(Y_n(s))}{\nabla_{\mathcal{V}_n} u_{\lambda,n}(Y_n(s))}_U\right]ds, \qquad t>0,
\end{align*}
and so
\begin{align*}
&\left|\int_0^te^{-\lambda s}\mathbb{E}\left[g(X_n(s))\right]ds-\int_0^te^{-\lambda s}\mathbb{E}\left[g(Y_n(s))\right]ds\right|\leq e^{-\lambda t}\mathbb{E}\left[\left|u_{\lambda,n}(X_n(t))-u_{\lambda,n}(Y_n(t))\right|\right]\\
&\phantom{aaaaaaaaaaa}+\int_0^te^{-\lambda s}\mathbb{E}\left[\left|\scal{B_n(X_n(s))-B_n(X(s))}{\nabla_{\mathcal{V}_n} u_{\lambda,n}(X_n(s))}\right|\right]ds \\
&\phantom{aaaaaaaaaaa}+\int_0^te^{-\lambda s}\mathbb{E}\left[\left|\scal{B_n(Y_n(s))-B_n(Y(s))}{\nabla_{\mathcal{V}_n} u_{\lambda,n}(Y_n(s))}\right|\right]ds.
\end{align*}
By Proposition \ref{Sol-n-int}, for every $t>0$ we get
\begin{align*}
&\left|\int_0^te^{-\lambda s}\mathbb{E}\left[g(X_n(s))\right]ds-\int_0^te^{-\lambda s}\mathbb{E}\left[g(Y_n(s))\right]ds\right|\leq e^{-\lambda t}2\norm{g}_{\infty}\\
&\phantom{aaaaaaaaaaa}+\int_0^te^{-\lambda s}\norm{g}_{\infty}\left(\mathbb{E}\left[\norm{B_n(X_n(s))-B_n(X(s))}\right]+\mathbb{E}\left[\norm{B_n(Y_n(s))-B_n(Y(s))}\right]\right)ds.
\end{align*}
Letting $n$ tend to $\infty$, by the uniformly continuity of $B$, \eqref{supE} and the dominated convergence theorem we obtain
\begin{align*}
&\left|\int_0^te^{-\lambda s}\mathbb{E}\left[g(X(s))\right]ds-\int_0^te^{-\lambda s}\mathbb{E}\left[g(Y(s))\right]ds\right|\leq e^{-\lambda t}2\norm{g}_{\infty}, \qquad t>0.
\end{align*}
finally letting $t$ go to $\infty$ we get
\begin{align*}
&\left|\int_0^\infty e^{-\lambda s}\mathbb{E}\left[g(X(s))\right]ds-\int_0^\infty e^{-\lambda s}\mathbb{E}\left[g(Y(s))\right]ds\right|=0.
\end{align*}
We have so proved that for every $x\in H$ and $g\in UC_b(H)$
\begin{align*}
\int_0^\infty e^{-\lambda s}\mathbb{E}\left[g(X(s))\right]ds=\int_0^\infty e^{-\lambda s}\mathbb{E}\left[g(Y(s))\right]ds,
\end{align*}
Hence by the properties of the Laplace transform, the proof is complete.
\end{proof}

\begin{proof}[Proof of Theorem \ref{WeakLoc}]
Theorem \ref{WeakLoc} is a immediate consequence of Theorem \ref{Weak1} and the method exploited in \cite[Section 6]{Pri2021}.
The only difference is that we consider the second-order Kolmogorov operator associated with \eqref{SDE} given by
\[
\mathcal{N}\varphi(x):=\frac{1}{2}{\rm Trace}\left[GG^*\nabla^2\varphi(x)\right]+\scal{x}{A^*\nabla\varphi(x)}_H+\scal{B(x)}{\mathcal{V}^*\nabla\varphi(x)}_H,\quad \varphi\in\xi_A(H),
\]
where $\xi_A(H)$ is the span of smooth bounded cylindrical functions given by
\[
\xi_A(H):={\rm span}\{ {\mbox{real and imaginary parts of functions }} \varphi:x\in H\rightarrow e^{i\scal{x}{h}} \mbox{ with } 
 h\in \cup_{n\in\N}H_n\}.
\]
We emphasize that this situation is even simpler than the one in \cite[Section 6]{Pri2021}, since $\mathcal{V}$ is bounded.
\end{proof}

\section{Examples}

\subsection{Stochastic heat equation}\label{Heat-case}
Let $U=H=L^2([0,2\pi]^d)$ with $d\in\N$. Consider the SPDE
\begin{align}\label{eqFObeta}
\left\{
\begin{array}{ll}
\displaystyle  dX(t)=\Delta X(t)dt+(-\Delta)^{-\beta}B(X(t))dt+(-\Delta)^{-\gamma}dW(t), \quad t\in[0,T],  \vspace{1mm} \\
X(0)=x\in H,
\end{array}
\right.
\end{align}
where $B\in UC(H;H)$, $\beta,\gamma\geq 0$ and $\Delta$ is the realization of the Laplace operator with Dirichlet boundary conditions in $H$. We recall that $\Delta$ is self-adjoint and there exists an orthonormal basis $\{e_k:k\in\N\}$ of $H$ consisting in eigenvectors of $\Delta$. Hence, the spaces $\{H_n\}_{n\in\N}$ defined by
\[
H_0=\{0\},\quad H_n:={\rm span}\{e_1,...,e_n\},\quad n\in\N,
\]
satisfy Hypotheses \ref{hyp:uniqueness}. Moreover, for every $k\in\N$, we have
\begin{equation*}
\Delta e_k=-\lambda_ke_k, \quad \lambda_k\sim k^{2/d}.
\end{equation*}
By easy computations, for every $t>0$ we obtain that
\begin{align*}
& Q_t:=\int_0^te^{2s\Delta}(-\Delta)^{-2\gamma}ds=\frac12(-\Delta)^{-(1+2\gamma)}(\Id_H-e^{2t\Delta}).
\end{align*}

By the same arguments used in the proof of \cite[Proposition 5.14]{AddBig}, we deduce the following proposition.
\begin{prop}\label{conti}
$ $
  \begin{enumerate}
      \item If $\gamma>\frac{d}{4}-\frac{1}{2}$, then \eqref{tra-con} hold. 
      \item There exists a constant $c>0$ such that for every  $t>0$ we have
\begin{align}
&\norm{\Gamma_t}_{\mathcal{L}(H)}\leq \frac{c}{t^{\frac12+\gamma}}\label{heat-controllo},\\
&\norm{\Gamma_t(-\Delta)^{-\beta}}_{\mathcal{L}(H)}\leq \frac{c}{t^{\frac12+\gamma-\beta}}.\notag
\end{align}
\end{enumerate}
  \end{prop}

Combining Theorem \ref{WeakLoc} and Proposition \ref{conti}, we obtain the following result.

\begin{thm}\label{Heat-weak}
Assume that 
\[
\gamma>\frac{d}{4}-\frac{1}{2},\qquad \beta>\gamma-\frac{1}{2}.
\]
Then weak uniqueness hold for \eqref{eqFObeta}.
\end{thm}

\begin{rmk}\label{rmk-heat}
In  \cite{Cho-Gol1995} and \cite{Kun2013} the weak uniqueness is proved assuming
    \[
\int^t_0 \|\Gamma_s\|_{\mathcal{L}(H;H)}ds<\infty.
\]
Therefore, taking into account \eqref{heat-controllo}, the results in \cite{Cho-Gol1995} and \cite{Kun2013} are applicable to \eqref{eqFObeta} only in the cases $d=1,2,3$.
\end{rmk}

\subsection{Stochastic damped equations}
\label{sub:stoch_damp_eq}
Let $U$ be a separable Hilbert space and let $H:=U\times U$. We deal with a semilinear stochastic differential equation of the form
\begin{align}\label{damped_stoc_equation}
\left\{
\begin{array}{ll}
\displaystyle \frac{\partial^2y}{\partial t^2}(t)
= -\Lambda y(t)-\rho\Lambda^{\alpha}\left(\frac{\partial y}{\partial t}(t)\right)+\Lambda^{-\beta}B\left(y(t),\frac{\partial y}{\partial t}(t)\right)+\Lambda^{-\gamma}d{W}(t), & t\in[0,T], \vspace{1mm} \\
y(0)=y_0, & \vspace{1mm} \\
\displaystyle \frac{\partial y}{\partial t}(0)=y_1,
\end{array}
\right.
\end{align}
where $\alpha\in(0,1)$, $\rho,\beta,\gamma\geq0$ and $B\in UC(H;U)$. We consider a self-adjoint operator $\Lambda:D(\Lambda)\subseteq U\to U$ of positive type with compact resolvent, simple eigenvalues $(\mu_n)_{n\in\N}$, blowing up as $n$ tends to infinity, and corresponding (non normalized) eigenvectors $\{e_n:n\in\N\}$. 

The operators  $A:D(A)\subseteq H \to H$ and $G,\mathcal{V}:U\to H$ in \eqref{SDE} are respectively defined as
\begin{align}
\label{damped_def_op_A_G_Lambda}
& D(A):= \left\{\begin{pmatrix}
h_1 \\ h_2    
\end{pmatrix}:h_2\in D(\Lambda^{\frac12}), \ h_1+\rho\Lambda^{\alpha-\frac12}h_2\in D(\Lambda^{\frac12})\right\}, \notag \\ 
& A:=\begin{pmatrix}
0 & \Lambda^{\frac12} \\
-\Lambda^{\frac12} & -\rho \Lambda^{\alpha}
\end{pmatrix}, 
\quad  G=\begin{pmatrix}
0 \\ \Lambda^{-\gamma}   
\end{pmatrix}, \quad \mathcal{V}=\begin{pmatrix}
 0 \\
 \Lambda^{-\beta}
\end{pmatrix}.
\end{align}
further $x=\begin{pmatrix}
\Lambda^{\frac12}y_0 \\ y_1    
\end{pmatrix}\in H$ with $y_0\in D(\Lambda^{\frac12})$ and $y_1\in U$, $X=\begin{pmatrix}
X_1 \\ X_2    
\end{pmatrix}$, $y=\Lambda^{-\frac12}X_1$ and $\frac{\partial y}{\partial t}=X_2$.
In this situation, if $\rho^2\neq 4\mu_n^{1-2\alpha}$, then the operator $A$, defined as in \eqref{damped_def_op_A_G_Lambda}, has simple eigenvalues
\begin{align}
\label{autovalori}
\lambda_n^{\pm}=\frac{-\rho\mu_n^{\alpha}\pm\sqrt{\rho^2\mu_n^{2\alpha}-4\mu_n}}{2}, \qquad \lambda_n^++\lambda_n^-=-\rho\mu_n^\alpha, \qquad \lambda_n^+\lambda_n^-=\mu_n, \qquad n\in\N,    
\end{align}
with corresponding eigenvectors
\begin{align}
\label{autovettoriA}
\Phi_n^+=\left(\begin{matrix}\mu_n^{\frac12} e_n \\ \lambda_n^+e_n\end{matrix}\right), \qquad \Phi_n^-=\chi_n\left(\begin{matrix}\mu_n^{\frac12} e_n \\
\lambda_n^-e_n\end{matrix}\right), \qquad n\in\N.    
\end{align}
The operator $A^*:D(A^*)\subseteq H\to H$ is defined as
\begin{align*}
A^*:=\begin{pmatrix}
0 & -\Lambda^{\frac12} \\
\Lambda^{\frac12} & -\rho \Lambda^{\alpha}
\end{pmatrix}, \qquad D(A^*)=D(A)
\end{align*}
and has simple eigenvalues $\lambda_n^{\pm}$, $n\in\N$, with corresponding eigenvectors
\begin{align}
\label{autovettoriA^*}
\Psi_n^+=\left(\begin{matrix}-\mu_n^{\frac12} e_n \\ \lambda_n^+e_n\end{matrix}\right), \qquad \Psi_n^-=\chi_n\left(\begin{matrix}-\mu_n^{\frac12} e_n \\
\lambda_n^-e_n\end{matrix}\right), \qquad n\in\N.    
\end{align}
It is possible to write $H=H^++H^-$ (non-orthogonal, direct sum), where
\begin{align*}
H^+=\overline{\{\Phi_n^+:n\in\N\}}=\overline{\{\Psi_n^+:n\in\N\}}, \qquad H^-=\overline{\{\Phi_n^-:n\in\N\}}=\overline{\{\Psi_n^-:n\in\N\}}.    
\end{align*}
Hence, any $h\in H$ can be written in a unique way as $h=h^++h^-$ with $h^+\in H^+$ and $h^-\in H^-$. Further, the operators $A$ and $e^{tA}$, $t\geq0$, act on $D(A)$ and $H$, respectively, as
\begin{align}
Ah= & \sum_{n=1}^\infty\left(\lambda_n^+\langle h^+,\Phi_n^+\rangle_H\Phi_n^++\lambda_n^-\langle h^-,\Phi_n^-\rangle_H\Phi_n^-\right), && h\in D(A),\label{operatoreA} \\
e^{tA}h
= & \sum_{n=1}^\infty\left(e^{\lambda_n^+t}\langle h^+,\Phi_n^+\rangle_H\Phi_n^++e^{\lambda_n^-t}\langle h^-,\Phi_n^-\rangle_H\Phi_n^-\right), && h\in H, \ t\geq0.
\label{semiA}
\end{align}
For every $v\in U$, the element $\left(\begin{matrix}
0 \\ v
\end{matrix}\right)$ can be written as
\begin{align*}
\left(
\begin{matrix}
0 \\ v
\end{matrix}
\right)= \sum_{n=1}^\infty\left(b_n^+ v_n\Phi_n^++b_n^- v_n\Phi_n^-\right), \qquad v_n:=\frac{\langle v,e_n\rangle_U}{\|e_n\|_U}, \ n\in\N,
\end{align*}
for suitable coefficients $b_n^+,b_n^-$, $n\in\N$. It follows that, for every $v\in U$, we have
\begin{align*}
\mathcal{V}v
= & \left(\begin{matrix}
0 \\ \Lambda^{-\beta} v    
\end{matrix}\right)
= \sum_{n=1}^\infty\mu_n^{-\beta} v_n\left(b_n^+\Phi_n^++b_n^-\Phi_n^-\right).
\end{align*}
From the above construction and \cite[formulae (2.3.14)-(2.3.18)]{Tri2003} it follows that, if $\alpha\in\left(0,\frac12\right]$, then \begin{align}
& {\rm Re}(\lambda_n^{+}), {\rm Re}(\lambda_n^{-})\sim -\mu_n^\alpha, \quad |\lambda_n^{+}|,|\lambda_n^{-}|\sim \mu_n^{\frac12}, \quad  \,  \|e_n\|_{U}\sim \mu_n^{-\frac12}, \quad |\lambda_n^+-\lambda_n^-|\sim \mu_n^{\frac12}, \notag \\
& b_n^{\pm}\sim {\rm const}(b), \quad \chi_n\sim {\rm const}(\chi),
\label{damped_stime_coefficienti_avl}
\end{align}
definitively with respect to $n\in\N$. The case $\alpha\in\left[\frac12,1\right)$ is analogously treated. We stress that, in this case, again from \cite[formulae (2.3.14)-(2.3.18)]{Tri2003}, the asymptotic behaviour in \eqref{damped_stime_coefficienti_avl} is replaced by
\begin{align}
& {\rm Re}(\lambda_n^+)\sim -\mu_n^{1-\alpha}, \quad {\rm Re}(\lambda_n^-)\sim -\mu_n^\alpha, \quad |\lambda_n^+|\sim \mu_n^{1-\alpha}, \quad |\lambda_n^-|\sim\mu_n^\alpha, \quad  \|e_n\|_{U}\sim \mu_n^{-\frac12}, \notag \\
& |\lambda^-_n-\lambda_n^+|\sim \mu_n^{\alpha}, \quad b_n^-\sim {\rm const}(b), \quad \chi_n,b_n^+\sim \mu_n^{\frac12-\alpha},
\label{damped_stime_coefficienti_avl_2}
\end{align}
definitively with respect to $n\in\N$. 

The spaces $\{H_n\}_{n\in\N}$, defined as 
$H_n:={\rm span}\{\Phi_k^+,\Phi_k^-\; :\; k\in\{1,\ldots n\}\}
={\rm span}\{\Psi_k^+,\Psi_k^-\; :\; k\in\{1,\ldots n\}\}$, $n\in\N$, fulfill Hypotheses \ref{hyp:uniqueness}(i). Moreover by \eqref{semiA}, 
\eqref{damped_stime_coefficienti_avl} and \eqref{damped_stime_coefficienti_avl_2} the semigroup $\{e^{tA}\}_{t\geq 0}$ is compact.

We consider the stochastic convolution
\begin{align*}
W_A(t):=\int_0^t e^{(t-s)A}GdW(s), \qquad t\geq0. 
\end{align*}
The following proposition shows that Hypotheses \ref{hyp:standard} hold true.
\begin{prop}{\cite[Proposition 5.3]{AddBig}}
\label{contidamped}
Let $A$ and $G$ be as in \eqref{damped_def_op_A_G_Lambda}. Assume that  there exist $\delta>0$ and a positive constant $c$ such that for every $n\in\N$ we have $\mu_n=c n^\delta$ and $\delta>\frac{1}{2\gamma+\alpha}$. Then \eqref{tra-con} hold.
\end{prop}

Now we prove that \eqref{supercontron} hold. To this aim, for every $t>0$ we consider the control problem
\begin{align}
\label{control_problem_n}
\dot{Y}(\tau)=AY(\tau)+Gu(\tau), \ \tau\in(0,t], \qquad Y(0)=h\in H,    
\end{align}
where $u:[0,t]\to U$. 
We recall that, if we consider the mild solution
\begin{align}
\label{mild_sol_problem_n}
Y(\tau)=e^{\tau A}h+\int_0^\tau e^{(\tau-s)A}Gu(s)ds, \qquad \tau\in[0,t],    
\end{align}
to problem \eqref{control_problem_n}, then $e^{tA}(H)\subseteq Q_{t}^{\frac12}(H)$ if and only if the system \eqref{control_problem_n} is null-controllable, i.e., for every $u\in L^2([0,t];U)$ problem \eqref{control_problem_n} admits a unique mild solution given by  \eqref{mild_sol_problem_n} and for every $h\in H$ there exists a control $u\in L^2([0,t];U)$ such that the corresponding mild solution $Y$ satisfies $Y(t)=0$ (see \cite{Zab08}). In this case, we say that $u$ steers $h$ to $0$  at time $t$ and, for every $t>0$ and $h\in H$, we have
\begin{align*}
\|Q_{t}^{-\frac12}e^{tA}h\|_{H}=\mathcal E_{C}(t,h):=\inf\{\|u\|_{L^2(0,t;U)}:u \textrm{ steers $h$ to $0$ at time $t$}\}. 
\end{align*}
In particular, for every $v\in U$ we have
$\|Q_{t}^{-\frac12}e^{tA}\mathcal Vv\|_{H}=\mathcal E_{C}(t,\mathcal Vv)$.

\begin{prop}\label{controllo}
Let $T>0$ and $\alpha\in (0,1)$. There exists a constant $c_{\alpha,T}>0$ such that
\begin{align}
\mathcal E_{C}(t,\mathcal Vv)
\leq \begin{cases}
\displaystyle \frac{c_{\alpha,T}\|v\|_{U}}{t^{\frac12+\frac{\gamma-\beta}{\alpha}\vee0}}, &  \ \alpha\in\left(0,\frac12\right], \\[2mm]
\displaystyle \frac{c_{\alpha,T}\|v\|_{U}}{t^{\frac12+\frac{\gamma-\beta}{1-\alpha}\vee0}}, & \ \alpha\in\left[\frac12,1\right),
\end{cases}\qquad t\in(0,T], \ v\in U.
\end{align}
In particular, for every $t\in [0,T]$, we have
\begin{align}
\norm{\Gamma_t\mathcal{V}}_{\mathcal{L}(U;H)}
\leq \begin{cases}
\displaystyle \frac{c_{\alpha,T}}{t^{\frac12+\frac{\gamma-\beta}{\alpha}\vee0}}, &  \ \alpha\in\left(0,\frac12\right], \\[2mm]
\displaystyle \frac{c_{\alpha,T}}{t^{\frac12+\frac{\gamma-\beta}{1-\alpha}\vee0}}, &  \ \alpha\in\left[\frac12,1\right).
\end{cases}
\end{align}
\end{prop}
\begin{proof}
We adopt the same strategy of the proof of \cite[Theorems 5.4 \& 5.6]{AddBig}, \cite[Theorem 4.6]{AddBig2} and \cite[Appendix]{MasPri2017}, which are inspired by \cite[Chapter 1, Proposition 1.3]{Zab08}. Since the computations are similar, we skip some details and refer to the proof of the quoted theorems. \\
We introduce the matrix $[G|AG]$ given by
\begin{align*}
[G|AG]
= \left(\begin{matrix}
0 & \Lambda^{\frac12-\gamma} \\ \Lambda^{-\gamma} & -\rho \Lambda^{\alpha-\gamma}\end{matrix}\right)
= \Lambda^{-\gamma}\left(
\begin{matrix}
0 & \Lambda^{\frac12} \\
{\rm Id}_{U} & -\rho \Lambda^{\alpha}
\end{matrix}
\right).
\end{align*}
The matrix $[G|AG]$ is, at least formally, invertible and its inverse matrix is
\begin{align*}
K=\Lambda^{\gamma}\left( \begin{matrix}
 \rho\Lambda^{\alpha-\frac12} & {\rm Id}_{U} \\
 \Lambda^{-\frac12} & 0
\end{matrix}    
\right)
=\left(\begin{matrix}
K_1 \\ K_2
\end{matrix}
\right).
\end{align*}
Let $t>0$. We define the control $u:[0,t]\to U$ as
\begin{align}
\label{definizione_controllo}
u(\tau)=
\begin{cases}
K_1\psi_t(\tau)+K_2\psi_t'(\tau), & \tau\in(0,t), \\[1mm]
0, & \tau=0 \textrm{ or }\tau=t,
\end{cases}
\end{align}
where $h\in H$, $\psi_t(\tau)=-\varphi_t(\tau)e^{\tau A}h$ for every $\tau\in[0,t]$ and $\varphi_t(\tau)=\overline c_m\tau^{m}(t-\tau)$ for every $\tau\in[0,t]$ with $\overline c_m$ being a normalizing constant and $m\in\N$ satisfying $m-2(\gamma-\beta)/\alpha>-1$, if $\alpha\in\left(0,\frac12\right]$, or $m-2(\gamma-\beta)/(1-\alpha)>-1$, if $\alpha\in\left[\frac12,1\right)$. From definition, we get $|\varphi_t(\tau)|\sim \frac{\tau^{m}}{t^{m+1}}$ and $|\varphi_t'(\tau)|\sim \frac{\tau^{m-1}}{t^{m+1}}$. \\
Arguing as in \cite[Theorem 5.4]{AddBig}, we deduce that $u$ steers $h$ to $0$ at time $t$, .i.e., the function $Y_n$, defined in \eqref{mild_sol_problem_n}, with $h\in H$ and $u$ as in \eqref{definizione_controllo}, is a mild solution to \eqref{control_problem_n} and $Y(t)=0$.

Let us estimate the $L^2$-norm of $u$. In the following, $c$ is a positive constant which may vary line by line. From \eqref{autovalori} and \eqref{autovettoriA}, we get
\begin{align*}
K_1\psi_t(\tau)
= & -\varphi_t(\tau)\Lambda^\gamma[\rho \Lambda^{\alpha-\frac12}(e^{\tau A_n}h)_1+(e^{\tau A_n}h)_2] \\
= & -\varphi_t(\tau)
\sum_{k\in\N}\mu_k^\gamma[\langle h^+,\Phi_k^+\rangle_H e^{\lambda_k^+ \tau}(\rho\mu_k^{\alpha}+\lambda_k^+)e_k+\langle h^-,\Phi_k^-\rangle_H e^{\lambda_k^- \tau}\chi_k(\rho\mu_k^{\alpha}+\lambda_k^-)e_k  ] \\
= & \varphi_t(\tau)
\sum_{k\in\N}\mu_k^\gamma[\langle h^+,\Phi_k^+\rangle_H e^{\lambda_k^+ \tau}\lambda_k^-e_k+\langle h^-,\Phi_k^-\rangle_H e^{\lambda_k^- \tau}\chi_k\lambda_k^+e_k  ], \qquad \tau\in[0,t].
\end{align*}

We now need an explicit expression for $\langle h^+,\Phi_k^+\rangle_H$ and $\langle h^-,\Phi_k^-\rangle_H$, $k\in\N$. From \eqref{autovettoriA} we infer that
\begin{align*}
h
= & \left(\begin{matrix}
h_1 \\ h_2
\end{matrix}\right) 
= \left(\begin{matrix}
\sum_{k\in\N}(h_1)_k\frac{e_k}{\|e_k\|_U} \\
\sum_{k\in\N}(h_2)_k\frac{e_k}{\|e_k\|_U}
\end{matrix}
\right)=\sum_{k\in\N}(\langle h^+,\Phi_k^+\rangle_H\Phi_k^++\langle h^-,\Phi_k^-\rangle_H\Phi_k^-) \\
= & \left(
\begin{matrix}
\sum_{k\in\N}\mu_k^{\frac12}(\langle h^+,\Phi_k^+\rangle_H+\chi_k\langle h^-,\Phi_k^-\rangle_H)e_k \\
\sum_{k\in\N}(\lambda_k^+\langle h^+,\Phi_k^+\rangle_H+\chi_k\lambda_k^-\langle h^-,\Phi_k^-\rangle_H)e_k
\end{matrix}
\right),
\end{align*}
where $(h_j)_k=\langle h_j,e_k\rangle_U\|e_k\|_{U}^{-1}$, $j=1,2$, for every $k\in\N$. By comparing the corresponding components of $e_k$, $k\in\N$, we infer that, for every $k\in\N$,
\begin{align*}
(h_1)_k
= & \mu_k^{\frac12}(\langle h^+,\Phi_k^+\rangle_H+\chi_k\langle h^-,\Phi_k^-\rangle_H)\|e_k\|_U,  \qquad
(h_2)_k
=  (\lambda_k^+\langle h^+,\Phi_k^+\rangle_H+\chi_k\lambda_k^-\langle h^-,\Phi_k^-\rangle_H)\|e_k\|_U.
\end{align*}
It follows that, for every $k\in\N$,
\begin{align}
\label{expl_h+h-}
\langle h^+,\Phi_k^+\rangle_H
= & \frac{\lambda_k^-(h_1)_k-\mu_k^{\frac12}(h_2)_k}{\mu_k^{\frac12}(\lambda_k^--\lambda_k^+)\|e_k\|_U},\qquad
\langle h^-,\Phi_k^-\rangle_H
=  \frac{\lambda_k^+(h_1)_k-\mu_k^{\frac12}(h_2)_k}{\chi_k\mu_k^{\frac12}(\lambda_k^+-\lambda_k^-)\|e_k\|_U}.
\end{align}
Notice that we are interested in $h=\mathcal{V}v$ with $v\in U$, so $h_1=0$ and
\begin{align*}
h_2=\Lambda^{-\beta}v
=\sum_{k=1}^n\mu_k^{-\beta}v_k \frac{e_k}{\|e_k\|_U}.
\end{align*}
By replacing the expression of $h_2$ in \eqref{expl_h+h-}, it follows that
\begin{align}
\langle h^+,\Phi_k^+\rangle_{H}
= & -\frac{\mu_k^{\frac12-\beta}v_k}{\mu_k^{\frac12}(\lambda_k^--\lambda_k^+)\|e_k\|_U}=-\frac{\mu_k^{-\beta}v_k}{(\lambda_k^--\lambda_k^+)\|e_k\|_U}, \qquad k\in\N, \label{stima_h^+}\\
\langle h^-,\Phi_k^-\rangle_H
= & \frac{\mu_k^{\frac12-\beta}v_k}{\chi_k\mu_k^{\frac12}(\lambda_k^--\lambda_k^+)\|e_k\|_{U}}
= \frac{\mu_k^{-\beta}v_k}{\chi_k(\lambda_k^--\lambda_k^+)\|e_k\|_{U}}, \qquad k\in\N.
\label{stima_h^-}
\end{align}
Substituting in $K_1\psi_t$, we infer that
\begin{align*}
K_1\psi_t(\tau)
= & \varphi_t(\tau)\sum_{k\in\N}\frac{\mu_k^{\gamma-\beta}}{(\lambda_k^--\lambda_k^+)\|e_k\|_U}[-e^{\lambda_k^+\tau}\lambda_k^-e_k+ e^{\lambda_k^-\tau}\lambda_k^+e_k] v_k, \qquad \tau\in(0,t),
\end{align*}
from which it follows that
\begin{align*}
\|K_1\psi_t(\tau)\|_U^2
= & (\varphi_t(\tau))^2\sum_{k\in\N}\frac{\mu_k^{2(\gamma-\beta)}}{|\lambda_k^--\lambda_k^+|^2}|e^{\lambda_k^+\tau}\lambda_k^-+e^{\lambda_k^-\tau}\lambda_k^+|^2v_k^2, \quad \tau\in(0,t).
\end{align*}
The estimates change if $\alpha\in\left(0,\frac12\right]$ or $\alpha\in\left[\frac12,1\right)$. In the first case, from \eqref{damped_stime_coefficienti_avl} we get
\begin{align*}
\|K_1\psi_t(\tau)\|_U^2
\leq & c(\varphi_t(\tau))^2\sum_{k=1}^n(|\mu_k^{\gamma-\beta}e^{-\rho\mu_k^{\alpha}\tau}|^2+|\mu_k^{\gamma-\beta}e^{-\rho\mu_k\tau^{\alpha}}|^2)v_k^2
\leq  \tau^{\left(-2\frac{\gamma-\beta}{\alpha}\right)\wedge 0}(\varphi_t(\tau))^2\sum_{k=1}^nv_k^2
\end{align*}
for every $\tau\in(0,t)$,
while in the second one formulae \eqref{damped_stime_coefficienti_avl_2} give
\begin{align*}
\|K_1\psi_t(\tau)\|_U^2
\leq & c(\varphi_t(\tau))^2\sum_{k\in\N}(\mu_k^{2(\gamma-\beta)}|e^{\lambda_k^+\tau}|^2+\mu_k^{2(\gamma-\beta)+2-4\alpha}|e^{\lambda_k^-\tau}|^2)v_k^2 \\
\leq & c(\varphi_t(\tau))^2\sum_{k\in\N}(\mu_k^{-2(\gamma-\beta)}e^{-2\rho\mu_k^{1-\alpha}\tau}+\mu_k^{2(\gamma-\beta+1-2\alpha)}e^{-2\rho\mu_k^{\alpha}\tau})v_k^2 \\
\leq & c(\varphi_t(\tau))^2(\tau^{(-2\frac{\gamma-\beta}{1-\alpha})\wedge 0}+\tau^{(-2\frac{\gamma-\beta+1-2\alpha}{\alpha})\wedge0})\sum_{k\in\N}v_k^2, \qquad \tau\in(0,t).
\end{align*}
Let us notice that, since $\alpha\in\left[\frac12,1\right)$, it follows that $\frac{\gamma-\beta}{1-\alpha}\vee0\geq \frac{\gamma-\beta+1-2\alpha}{\alpha}\vee0$, due to the fact that, if $\gamma-\beta\geq0$, then the right-hand side is smaller than $\frac{\gamma-\beta}{\alpha}\leq \frac{\gamma-\beta}{1-\alpha}$, and if $\gamma-\beta<0$ then both the sides are $0$. Hence,
\begin{align}
\|K_1\psi_t(\tau)\|_U^2
\leq \begin{cases}
c(\varphi_t(\tau))^2\tau^{\left(-2\frac{\gamma-\beta}{1-\alpha}\right)\wedge 0}\|v\|_{U}^2, & \tau\in(0,t), \ \alpha\in\left[\frac12,1\right), \\[1mm]
c(\varphi_t(\tau))^2\tau^{\left(-2\frac{\gamma-\beta}{\alpha}\right)\wedge 0}\|v\|_{U}^2, & \tau\in(0,t), \  \alpha\in\left(0,\frac12\right].
\end{cases}
\label{stima_K_1}
\end{align}
\\
Now we estimate the $U$-norm of $K_2\psi_t'(\tau)$. To this aim, we notice that, from \eqref{autovalori} and \eqref{autovettoriA},
\begin{align*}
K_2\psi'_t(\tau)
= & -\Lambda^{\gamma-\frac12}_n(\varphi_t'(\tau)(e^{\tau A_n}h)_1+\varphi_t(\tau)(A_ne^{\tau A_n}h)_1) \\
= & -\sum_{k\in\N} \mu_k^{\gamma}[\varphi_t'(\tau)(\langle h^+,\Phi_k^+\rangle e^{\lambda_k^+\tau}e_k+\chi_k\langle h^-,\Phi_k^-\rangle e^{\lambda_k^-\tau}e_k) \\
& + \varphi_t(\tau)(\langle h^+,\Phi_k^+\rangle e^{\lambda_k^+\tau}\lambda_k^+e_k+\chi_k\langle h^-,\Phi_k^-\rangle e^{\lambda_k^-\tau}\lambda_k^-e_k)], \qquad \tau\in(0,t).
\end{align*}
Taking \eqref{stima_h^+} and \eqref{stima_h^-} into account, we infer that
\begin{align*}
\|K_2\psi_t'(\tau)\|_U^2
= & \sum_{k\in\N}\mu_k^{2(\gamma-\beta)}\frac{|\varphi_t'(\tau)(-e^{\lambda_k^+\tau}+e^{\lambda_k^-\tau})+\varphi_t(\tau)(-\lambda_k^+e^{\lambda_k^+\tau}+\lambda_k^-e^{\lambda_k^-\tau})|^2}{|\lambda_k^--\lambda_k^+|^2}v_k^2, \quad \tau\in(0,t).
\end{align*}
As above, we split the cases $\alpha\in\left(0,\frac12\right]$ and $\alpha\in\left[\frac12,1\right)$. If $\alpha\in\left(0,\frac12\right]$, then from \eqref{damped_stime_coefficienti_avl} we get
\begin{align}
& |\varphi_t'(\tau)|^2\left|\frac{-e^{\lambda_k^+\tau}+e^{\lambda_k^-\tau}}{\lambda_k^--\lambda_k^+}\right|^2
\leq c\frac{\tau^{2(m-1)}|e^{\lambda_k^+\tau}|^2}{t^{2(m+1)}}\left|\frac{e^{-\lambda_k^+\tau+\lambda_k^-\tau}-1}{\lambda_k^--\lambda_k^+}\right|^2
\leq c\frac{\tau^{2m}}{t^{2(m+1)}}e^{-2\rho\mu_k^{\alpha}}, && k\in\N, \label{forma_K_2_1}\\
& |\varphi_t(\tau)|^2\left|\frac{-\lambda_k^+e^{\lambda_k^+\tau}+\lambda_k^-e^{\lambda_k^-\tau}}{\lambda_k^--\lambda_k^+}\right|^2
\leq c\frac{\tau^{2m}}{t^{2(m+1)}}
\frac{|\lambda_k^+e^{\lambda_k^+\tau}|^2+|\lambda_k^-e^{\lambda_k^-\tau}|^2}{|\lambda_k^--\lambda_k^+|^2}
\leq c\frac{\tau^{2m}}{t^{2(m+1)}}e^{-2\rho\mu_k^\alpha}, && k\in\N
\label{forma_K_2_2}
\end{align}
for every $\tau\in(0,t)$, where we have used the fact that ${\rm Re}(\lambda_k^--\lambda_k^+)\leq 0$ definitively, which implies that
\begin{align*}
\|K_2\psi_t'(\tau)\|_{U}^2
\leq & c\sum_{k\in\N}\frac{\tau^{2m}}{t^{2(m+1)}}\mu_k^{2(\gamma-\beta)}e^{-2\rho\mu_k^{\alpha}}v_k^2 
\leq c\sum_{k\in\N}\frac{\tau^{2m+\left(-\frac{2(\gamma-\beta)}{\alpha}\right)\wedge 0}}{t^{2(m+1)}}v_k^2, \qquad \tau\in(0,t).
\end{align*}
If $\alpha\in\left[\frac12,1\right)$ and arguing as for \eqref{forma_K_2_1} and \eqref{forma_K_2_2}, then from \eqref{damped_stime_coefficienti_avl_2}, \eqref{stima_h^+} and \eqref{stima_h^-} it follows that
\begin{align*}
& |\varphi_t'(\tau)|^2\left|\frac{-e^{\lambda_k^+\tau}+e^{\lambda_k^-\tau}}{\lambda_k^--\lambda_k^+}\right|^2
\leq c\frac{\tau^{2m}}{t^{2(m+1)}}e^{-2\rho\mu_k^{1-\alpha}\tau}, && k\in\N, \\
& |\varphi_t(\tau)|^2\left|\frac{-\lambda_k^+e^{\lambda_k^+\tau}+\lambda_k^-e^{\lambda_k^-\tau}}{\lambda_k^--\lambda_k^+}\right|^2
\leq c\frac{\tau^{2m}}{t^{2(m+1)}}(\mu_k^{2(1-2\alpha)}e^{-2\rho\mu_k^{1-\alpha}\tau}+e^{-2\rho\mu_k^{\alpha}\tau}), && k\in\N
\end{align*}
for every $\tau\in(0,t)$, which gives
\begin{align*}
\|K_2\psi_t'(\tau)\|^2_U
\leq & c\frac{{\tau}^{2m}}{t^{2(m+1)}}(\tau^{\left(-\frac{2(\gamma-\beta)}{1-\alpha}\right)\wedge 0}+\tau^{\left(-\frac{2(\gamma-\beta+1-2\alpha)}{1-\alpha}\right)\wedge0}+\tau^{\left(-\frac{2(\gamma-\beta)}{\alpha}\right)\wedge 0})\sum_{k\in\N}v_k^2 \\
\leq & c\frac{\tau^{2m+\left(-\frac{2(\gamma-\beta)}{1-\alpha}\right)\wedge 0}}{t^{2(m+1)}}\sum_{k\in\N}v_k^2, \qquad \tau\in(0,t),
\end{align*}
where we have used the fact that $\frac{2(\gamma-\beta)}{\alpha}\vee0\leq \frac{2(\gamma-\beta)}{1-\alpha}\vee0$ and $ \frac{2(\gamma-\beta+1-2\alpha)}{1-\alpha}\leq \frac{2(\gamma-\beta)}{1-\alpha}$ for every $\alpha\in\left[\frac12,1\right)$.
Hence,
\begin{align}
\label{stima_K_2}
\|K_2\psi_t'(\tau)\|_U^2
\leq c\frac{t^{2m}}{t^{2m+2}}
\begin{cases}
\tau^{\left(-2\frac{\gamma-\beta}{1-\alpha}\right)\wedge 0}\|v\|_{U}^2, & \tau\in(0,t), \ \alpha\in\left[\frac12,1\right), \\[1mm]
\tau^{\left(-2\frac{\gamma-\beta}{\alpha}\right)\wedge 0}\|v\|_{U}^2, & \tau\in(0,t), \  \alpha\in\left(0,\frac12\right].
\end{cases}
\end{align}
If $\alpha\in\left(0,\frac12\right]$, then from \eqref{stima_K_1} and recalling that $|\varphi_t(\tau)|\sim \frac{\tau^{m}}{t^{m+1}}$, $\tau\in(0,t)$, it follows that
\begin{align*}
\|K_1\psi_t\|_{L^2(0,t;U)}^2
\leq & \frac{c}{t^{2m+2}}\int_0^t\tau^{2m+\left(-\frac{2(\gamma+\beta)}{\alpha}\right)\wedge 0}d\tau\|v\|_{U_n}^2
\leq ct^{\left(-\frac{2(\gamma+\beta)}{\alpha}\right)\wedge 0-1}\|v\|_{U_n}^2,
\end{align*}
and, taking \eqref{stima_K_2} into account, 
\begin{align*}
\|K_2\psi'_t\|_{L^2(0,t;U)}^2
\leq ct^{\left(-\frac{2(\gamma-\beta)}{\alpha}\right)\wedge 0-1}\|v\|_{U}^2.    
\end{align*}
If $\alpha\in\left[\frac12,1\right)$, then arguing as above we get
\begin{align*}
\|K_1\psi_t\|_{L^2([0,t];U)}^2
\leq ct^{\left(-\frac{2(\gamma-\beta)}{1-\alpha}\right)\wedge 0-1}\|v\|_{U}^2, \qquad 
\|K_2\psi'_t\|_{L^2([0,t];U)}^2
\leq ct^{\left(-\frac{2(\gamma-\beta)}{1-\alpha}\right)\wedge 0-1}\|v\|_{U}^2,    
\end{align*}
which gives the thesis.
\end{proof}

Combining Theorem \ref{WeakLoc} and Propositions \ref{contidamped} and \ref{controllo}, we obtain the following result.

\begin{thm}
\label{thm:damped_main_result_alpha<12}
Assume that:
\begin{enumerate}[\rm(i)]
\item $\alpha\in (0,1)$ and $B\in UC(H;U)$;
\item there exist $\delta>0$ and a positive constant $c$ such that for every $n\in\N$ we have $\mu_n=c n^\delta$ and $\gamma>\frac{1-\alpha\delta}{2\delta}$;
\item $\beta\geq 0$ verifies
\[
\begin{cases}
\displaystyle \beta>\gamma-\frac{\alpha}{2},\qquad\ &{\rm if}\ \alpha\in\left(0,\frac12\right], \\[2mm]
\displaystyle \beta>\gamma-\frac{1-\alpha}{2},\qquad\ &{\rm if}\ \ \alpha\in\left[\frac12,1\right).
\end{cases}
\]
\end{enumerate}
Then weak uniqueness hold for \eqref{damped_def_op_A_G_Lambda}.
\end{thm}
Let $d\in\N$ and let $\Delta$ be the realization in $L^2([0,2\pi]^d)$ of the Laplace operator with Dirichlet boundary conditions. The two next results follow by applying Theorem \ref{thm:damped_main_result_alpha<12} in two specific cases.

\begin{coro}[Stochastic damped wave equation]\label{Wave}
Let $d\in\N$. Assume that $U=L^2([0,1]^d)$, $\Lambda=-\Delta$, $B\in UC(H;U)$, $\alpha\in (0,1)$, $\beta,\gamma\geq 0$ and 
\begin{align*}
  \gamma> \frac{d-2\alpha}{4},\qquad  \begin{cases}
\displaystyle \beta>\gamma-\frac{\alpha}{2},\quad\ &{\rm if}\ \alpha\in\left(0,\frac12\right], \\[2mm]
\displaystyle \beta>\gamma-\frac{1-\alpha}{2},\quad\ &{\rm if}\ \ \alpha\in\left[\frac12,1\right).
\end{cases}
\end{align*}
Then weak uniqueness hold for \eqref{damped_def_op_A_G_Lambda}.
\end{coro}

\begin{coro}[Stochastic damped Euler-Bernoulli beam equation]\label{EB}
Let $d\in\N$. Assume that $U=L^2([0,1]^d)$, $\Lambda=\Delta^2$, $B\in UC(H;U)$, $\alpha\in (0,1)$, $\beta,\gamma\geq 0$ and 
\begin{align*}
  \gamma> \frac{d-4\alpha}{8},\qquad  \begin{cases}
\displaystyle \beta>\gamma-\frac{\alpha}{2},\quad\ &{\rm if}\ \alpha\in\left(0,\frac12\right], \\[2mm]
\displaystyle \beta>\gamma-\frac{1-\alpha}{2},\quad\ &{\rm if}\ \ \alpha\in\left[\frac12,1\right).
\end{cases}
\end{align*}
Then weak uniqueness hold for \eqref{damped_def_op_A_G_Lambda}.
\end{coro}

\begin{rmk}\label{rmk-damped}
We note that in both the cases of Corollary \ref{Wave} and Corollary \ref{EB}, by \cite[Proposition 5.4]{AddBig}, for every $T>0$ there exists a constant $\overline c>0$ such that, for every $t\in(0,T]$,
\[
\norm{\Gamma_t}_{\mathcal{L}(H;H)}
\leq 
\begin{cases}
\displaystyle \frac{\overline c}{t^{\frac12+(\gamma+\alpha-\frac12)/(1-\alpha)}}, & \gamma+2\alpha\geq \frac32, \\[4mm] 
\displaystyle \frac{\overline c}{t^{\frac32}}, & \gamma+2\alpha< \frac32. 
\end{cases}
\]
Hence, the results in \cite{Cho-Gol1995,Kun2013} are not applicable, since \eqref{K-contro} is not verified.
\end{rmk}

\end{document}